\documentclass{article}
\usepackage{amssymb}
\usepackage{graphicx}
\usepackage{amsmath}

\newtheorem{theorem}{Theorem}

\newtheorem{corollary}[theorem]{Corollary}

\newtheorem{lemma}[theorem]{Lemma}

\newtheorem{proposition}[theorem]{Proposition}

\renewcommand{\a}{{\bf{a}}}

\begin{document}

\title{Generic groups acting on regular trees}
\author{Mikl\'{o}s Ab\'{e}rt and Yair Glasner \thanks{Some of the work was
carried out at the Institute for Advanced Studies and supported by NSF grant
DMS-0111298}}
\maketitle

\begin{abstract}
Let $T$ be a $k$-regular tree ($k\geq 3$) and $A=\mathrm{Aut}(T)$ its
automorphism group. We analyze a generic finitely generated subgroup $\Gamma
$ of $A$. We show that $\Gamma $ is free and establish a trichotomy on the
closure $\overline{\Gamma }$ of $\Gamma $ in $A.$ It turns out that $%
\overline{\Gamma }$ is either discrete, compact or has index at most $2$ in $%
A$.
\end{abstract}

\section{Introduction}

Let $T$ be a $k$-regular tree ($k\geq 3$) and let $\mathrm{Aut}(T)$ be its
automorphism group. The topology of pointwise convergence turns $\mathrm{Aut}%
(T)$ into a locally compact, totally disconnected, unimodular topological
group. Let $\mu $ be a Haar measure on $\mathrm{Aut}(T)$ normalized so that
vertex stabilizers have measure $1$. Since $\mathrm{Aut}(T)$ is not compact,
$\mu $ is an infinite measure. The group $\mathrm{Aut}(T)$ has a unique
normal subgroup $\mathrm{Aut}^{0}(T)$ of index $2$, which is simple.

Let $n\geq 2$ be an integer. For the $n$-tuple $\a=(a_{1},a_{2},\ldots
,a_{n})\in \mathrm{Aut}(T)^{n}$ let $\langle \a\rangle =\langle
a_{1},a_{2},\ldots ,a_{n}\rangle $ denote the group generated by the $a_{i}$%
. We say that a group theoretic property $P$ is \textit{measure generic} if $%
P$ holds for the group $\langle \a\rangle $ for $\mu ^{n}$-almost all $\a\in
\mathrm{Aut}(T)^{n}$. Similarly, $P$ is said to be \textit{topologically
generic} if it holds for the group $\langle \a\rangle $ for all but a meager
(or first category) subset of $\mathrm{Aut}(T)^{n}$. The property $P$ is
called \textit{generic} if it is both measure and topologically generic.
E.g. the sentence ``A generic finitely generated subgroup of $\mathrm{Aut}(T)
$ is infinite'' means that for all $n\in \mathbb{N}$ the set $\{\a\in
\mathrm{Aut}(T)^{n}|\langle \a\rangle \text{ is finite}\}$ is a meager
nullset.

The main result of this paper is the following.

\begin{theorem}
\label{thm:main} Let $n\geq 2$ and let $\Gamma <\mathrm{Aut}(T)$ be a
generic subgroup on $n$ generators. Then $\Gamma $ is isomorphic to $F_{n}$,
the free group of rank $n$ and it falls into exactly one of the following
categories: \newline
a) $\Gamma $ is discrete;\newline
b) $\Gamma $ is precompact, that is, $\Gamma $ fixes a vertex or a geometric
edge of $T$;\newline
c) $\Gamma $ is dense in $\mathrm{Aut}(T)$ or in $\mathrm{Aut}^{0}(T)$.
\end{theorem}

We also show that all of a), b) and c) hold on a subset of infinite measure.

An essential tool we use in proving Theorem \ref{thm:main} is the natural
action of the automorphism group $\mathrm{Aut}(F_{n})$ on $\mathrm{Aut}%
(T)^{n}$ by measure-preserving homeomorphisms. An important tool is a
theorem of Weidmann implying that every non-discrete finitely generated free
subgroup of $\mathrm{Aut}(T)$ has a primitive element that fixes a vertex or
a geometric edge.

The group $\mathrm{Aut}(T)$ and its discrete subgroups have been extensively
studied in the literature (see \cite{Tits:Arbre}, \cite
{FtN:Harmonic_analysis_trees}, \cite{serre}, \cite{basslub} and references
therein). Precompact subgroups, which are just groups acting on rooted
trees, have also been much investigated (see \cite{grigor} and for random
generation \cite{bhatta}, \cite{AV-dimension_theory} and references therein).
The existence of finitely generated dense free subgroups in $\mathrm{Aut}(T)$
is a new phenomenon however, and it leads to a new example in the realm of infinite
permutation groups.

\begin{corollary}
\label{szep}For all $n\geq 2$ the free group $F_{n}$ has a primitive, but
not $2$-transitive action on a countable set $X$, such that for each finite
subset $Y\subseteq X$ the pointwise stabilizer $F_{Y}$ is nontrivial and
every nontrivial subnormal subgroup of $F$ is transitive on $X$.
\end{corollary}

It is natural to ask whether the trichotomy established in Theorem \ref
{thm:main} works for arbitrary locally compact groups. Namely, if is it true
that in such a group $G$, the closure of a generic subgroup is either open,
compact or discrete. For connected semisimple Lie groups (like $\mathrm{SL}%
(2,\mathbb{R})$) there is actually a dichotomy: generically, the closure is
either open or discrete (see \cite{gelander} for a proof in the compact case).
However, as we show in Proposition \ref{bzzz}, the trichotomy is not true
in general.

\bigskip

The paper is organized as follows. Section \ref{generalsec} introduces some
notation and known results that will be needed later. In Section \ref{densesec}
we show that a fixed hyperbolic and a generic elliptic element generate a
dense subgroup. This is the core result leading to Theorem \ref{thm:main}
that is proved in Section \ref{trichotomysec}. Finally, in Section \ref
{primitivsec} we apply our results, in particular, we obtain Corollary \ref
{szep} and Proposition \ref{bzzz}.

\bigskip

\noindent \textbf{Acknowledgement.} Parts of the paper have been discussed
at the Secret Seminar, for which we are openly grateful.

\section{Tree automorphisms \label{generalsec}}

Let $k\geq 3$ and let $T=T_{k}$ be a $k$-regular tree. For $x,y\in T$ let $%
[x,y]$ be the unique simple path going from $x$ to $y\ $and let $d(x,y)$ be
the length of $[x,y]$, the distance of $x$ and $y$ in the graph metric. If $%
x $ and $y$ are neighbours, let
\begin{equation*}
\mathrm{Shadow}_{x\longrightarrow y}=\left\{ t\in T\mid y\in \lbrack
x,t]\right\}
\end{equation*}
be the \emph{shadow of the directed edge }$x\longrightarrow y$. For $t\in T$
and $n\in N$ let $B(t,n)=\left\{ x\in T\mid d(t,x)\leq n\right\} $.

Let $A=\mathrm{Aut}(T)$. Taking stabilizers of finite subsets of $T$ as a
base of neighbourhoods of $1$, $A$ turns into a totally disconnected,
locally compact, unimodular topological group. Thus $A$ admits a left and
right invariant Haar measure $\mu $ and by normalizing we can assume that
every vertex stabilizer has $\mu $-measure $1$.

The tree $T$ is a bipartite graph and the stabilizer of the bipartite
partition, $A_{0}=\mathrm{Aut}^{0}(T)$ is the unique nontrivial normal
subgroup of $A$.

It is easy to see that a subgroup of $A$ is discrete if and only if it has
finite vertex stabilizers. Also, a subgroup is precompact, i.e., it has a
compact closure in $A$, if and only if it fixes a vertex or a geometric
edge. By a geometric edge we mean a set of two adjacent vertices.

Let us briefly recall the classification of tree automorphisms. Elements of $%
A$ are either \emph{elliptic} (fixing a vertex of $T$), \emph{inversions}
(fixing a geometric edge but not a vertex of $T$) or \emph{hyperbolic}
elements (the rest). For convenience, we want to treat inversions as
elliptic elements. To do so, let us take the so-called barycentric
subdivision of the tree $T$, that is, let us insert a new vertex in the
center of each edge of $T$. The new tree obtained this way is bi-regular and
its automorphism group equals $A$. Indeed, the only thing that changes is
that now inversions also fix a vertex, that is, they are elliptic. So, from
now on, $T$ denotes this geometric realization, and $A$ its automorphism
group. The price we pay for this convenience is that we need to prove all
our results in the realm of bi-regular trees.

For $a\in A$ let us define the \emph{minimal translation length} as
\begin{equation*}
l(a)=\min_{x\in T}\{d(x,xa)\}
\end{equation*}
and the \emph{support} as
\begin{equation*}
X(a)=\{x\in T|d(x,xa)=l(a)\}.
\end{equation*}
So, $a$ is elliptic if and only if $l(a)=0$. In this case $X(a)$ is a convex
subtree of $T$ consisting of the fixed points of $a$. If $a$ is hyperbolic,
then $X(a)$ is an $a$-invariant infinite geodesics called the \emph{axis} of
$a$ and $a$ acts on $X(a)$ by translation of length $l(a)$.

Now we will quote some classical results on groups acting on trees. The
reader can find them e.g. in Serre's book on trees \cite{serre}.

The set of elliptic elements is clearly open and therefore so is the set of
hyperbolics. It is also well understood which subgroups consist of only
elliptic or hyperbolic elements. A group consists of elliptic elements if
and only if it fixes a vertex. This follows from the fact that if $a,b$ are
elliptic with no common fixed point then $ab$ is hyperbolic.

Let $X$ and $Y$ be convex subsets of $T$. Let us define the \emph{projection
of }$Y$\emph{\ onto }$X$ to be the set
\begin{equation*}
\mathrm{Proj}_{X}Y=\left\{ x\in X\mid d(x,Y)=d(X,Y)\right\} .
\end{equation*}

Let $H=\{h_{1},h_{2},\ldots ,h_{n}\}$ be a set of hyperbolic elements and
let $X_{i}=X(h_{i})$. Then $H$ satisfies the \emph{Schottky condition} if
for every $1\leq i\leq n$ there exists a geodesic $Y_{i}\subset X_{i}$ of
length at most $l(h_{i})-1$ such that
\begin{equation*}
\bigcup_{j\neq i}\mathrm{Proj}_{X_{i}}(X_{j})\subseteq Y_{i}
\end{equation*}

Now a theorem of Lubotzky \cite{Lubotzky:LatticesOverLf} asserts that a
finitely generated subgroup of $A$ is discrete and free if and only if it
has a set of free generators satisfying the Schottky condition. This shows
that a subgroup of $A$ consists only of hyperbolic elements if and only if
it is discrete and free and that for a finite set of generators this is an
open condition.

\bigskip

Now we will fix our convention to treat automorphisms of rooted trees. Let $%
d_{0},d_{1},\ldots ,d_{n}$ be a sequence of positive integers ($n$ might be
infinite). Let $Z_{r}$ denote the set of residue classes modulo $r$. For $%
l<n+1$ let
\begin{equation*}
T_{l}=\left\{ (s_{0},s_{1},\ldots ,s_{l-1})\mid s_{i}\in Z_{d_{i}}\right\}
\end{equation*}
The \emph{root} is the empty sequence $()\in T_{0}$. Let
\begin{equation*}
T=T(d_{0},d_{1},\ldots ,d_{n})=\bigcup_{0\leq l<n+1}T_{l}
\end{equation*}
The level of a vertex $v\in T_{l}$ is defined to be $\left| v\right| =l$. A
vertex $w=(r_{0},r_{1},\ldots ,r_{m})\in T$ is a \emph{descendant }of $%
v=(s_{0},s_{1},\ldots ,s_{l})\in T$ (or $w>v$) if $l<m$ and $r_{i}=s_{i}$ ($%
0\leq i\leq l$). The vertex $w$ is a \emph{child} of $v$ if it is a
descendant and $m=l+1$. Drawing edges along the child relation turns $T$
into a spherically homogeneus tree rooted at $()$.

Let $v\in T$ be a vertex of level $l$. An \emph{elementary automorphism at }$%
v$ is defined as a permutation $\alpha \in Sym(Z_{d_{l}})$, the symmetric
group on $Z_{d_{l}}$ acting on $T$ as follows. For $w=(r_{0},r_{1},\ldots
,r_{m})\in T$ let
\begin{equation*}
w^{\alpha }=\left\{
\begin{array}{cc}
(r_{0},\ldots ,r_{l-1},r_{l}^{\alpha },r_{l+1},\ldots ,r_{m}) & \text{if }%
w\geq v \\
w & \text{otherwise}
\end{array}
\right.
\end{equation*}
Elementary automorphisms generate the full automorphism group $\mathrm{Aut}%
(T)$ (if $n$ is infinite, then one has to consider topological generation).
Furthermore, every automorphism $\varphi \in \mathrm{Aut}(T)$ can be
uniquely written as a product of elementary automorphisms
\begin{equation}
\varphi =\prod_{v\in T}\varphi (v)  \label{felbontas}
\end{equation}
where $\varphi (v)$ is an elementary automorphism at $v$ and the product is
taken in non-decreasing order with respect to the level of $v$. If $n$ is
infinite then the above infinite product should be understood as a limit.

Note that elementary automorphisms at incomparable vertices commute, so we
can rearrange our product accordingly. We will choose a rearrangement that
reflects looking at $\mathrm{Aut}(T)$ as the wreath product of $\mathrm{Aut}%
(T(d_{l},d_{l+1},\ldots ,d_{n}))$ with $\mathrm{Aut}(T(d_{0},d_{1},\ldots
,d_{l-1}))$ as follows. For $1\leq l<n$ and $\varphi \in \mathrm{Aut}(T)$
let
\begin{equation*}
\varphi _{l}=\prod_{\left| v\right| <l}\varphi (v)
\end{equation*}
and for $w\in T$ of level $l$ let
\begin{equation*}
\varphi _{w}=\prod_{v\geq w}\varphi (v)
\end{equation*}
again in non-decreasing order with respect to the level of $v$. Then
\begin{equation*}
\varphi =\varphi _{l}\prod_{\left| v\right| =l}\varphi _{v}
\end{equation*}
where the product is now in arbitrary order since the different $\varphi
_{v} $-s commute. The autmorphism $\varphi _{l}$ acts on the union of the
first $l $ levels of $T$ so it can be looked at as an element of $\mathrm{Aut%
}(T(d_{0},d_{1},\ldots ,d_{l-1}))$. The automorphism $\varphi _{v}$ acts on
the subtree $T(d_{l},d_{l+1},\ldots ,d_{n})$; this action is called the
\emph{state of }$\varphi $\emph{\ at }$v$\emph{\ }and is also denoted by $%
\varphi _{v}$.

\bigskip

Let $G$ be any group and let $n$ be an integer. Then $\mathrm{Aut}(F_{n})$
has a natural action on $\mathrm{Hom}(F_{n},G)$ as follows. For $\varphi \in
\mathrm{Aut}(F_{n})$ and $f\in \mathrm{Hom}(F_{n},G)$ let
\begin{equation*}
f^{\varphi }:w\longmapsto f(w^{\varphi })\text{ (}w\in F_{n}\text{)}
\end{equation*}
By fixing a minimal generating set for $F_{n}$ we can look at $\mathrm{Hom}%
(F_{n},G)$ as the set of $n$-tuples from $G$. Hence $\mathrm{Aut}(F_{n})$
acts on $G^{n}$ as well. This action can be best understood by Nielsen
transformations, the action of the Nielsen generators of $\mathrm{Aut}(F_{n})
$. Let $\a=(a_{1},a_{2},\ldots ,a_{n})\in G^{n}$. Then a \emph{right Nielsen
transformation of }$G^{n}$ is of the form
\begin{equation*}
(a_{1},\ldots ,a_{n})^{R_{i,j}^{\pm }}=(a_{1},\ldots
,a_{i-1},a_{i}a_{j}^{\pm 1},a_{i+1},\ldots ,a_{n})\text{ (}i\neq j\text{)}
\end{equation*}
and a \emph{left Nielsen transformation of }$G^{n}$ is
\begin{equation*}
(a_{1},\ldots ,a_{n})^{L_{i,j}^{\pm }}=(a_{1},\ldots ,a_{i-1},a_{j}^{\pm
1}a_{i},a_{i+1},\ldots ,a_{n})\text{ (}i\neq j\text{)}
\end{equation*}
Lastly, permutations of coordinates are also Nielsen transformations.

It turns out that if the generators $\a$ are chosen uniform randomly, then
the image will still be distributed in the same way. This is true in both
the measure theoretic and the topological settings by the following.

\begin{lemma}
\label{lem:Nielsen}Let $G$ be a unimodular locally compact topological group
with a Haar measure $\mu $. Then Nielsen transformations act on $G^{n}$ by $%
\mu $-preserving homeomorphisms.
\end{lemma}

\noindent \textbf{Proof.} For permutations of coordinates, the lemma is
trivial. Clearly $R_{i,j}^{\pm }$ and $L_{i,j}^{\pm }$ are continuous
bijections of $G^{n}$. Since their inverses are also Nielsen
transformations, they are homeomorphisms of $G^{n}$. To demonstrate that
Nielsen transformations are measure preserving it is enough to check the $2$
variable case $N:(x,y)\longmapsto (xy,y)$. Using Fubini's theorem, we have
\begin{eqnarray*}
\int_{G\times G}f\circ N(x,y)d\mu ^{2}(x,y) &=&\int_{G}\int_{G}f(xy,y)d\mu
(x)d\mu (y) \\
&=&\int_{G}\int_{G}f(x,y)d\mu (x)d\mu (y) \\
&=&\int_{G\times G}f(x,y)d\mu ^{2}(x,y)
\end{eqnarray*}
Note that we only used the right invariance of $\mu $ for this calculation.
For left Nielsen transformations we have to use the left invariance of $\mu $%
. $\square $

\bigskip

We say that two $n$-tuples $a,b\in G^{n}$ are \emph{Nielsen-equivalent} if
there is a sequence of Nielsen transformations leading from $a$ to $b$, or,
equivalently, if they lie in the same $\mathrm{Aut}(F_{n})$-orbit.

A crucial fact about Nielsen transformations is that they do not change the
subgroup the $n$-tuple generates. On the other hand, they do change the
tuple considerably, allowing us to find the `right' generating set for the
subgroup. This method is called the Nielsen method; it was first used by
Nielsen to show that finitely generated subgroups of a free group are free
\cite{Nielsne:Subgroups}. Since then, the method has been further developed
by Zieschang, Weidmann, Kapovich and others. We will need the following
theorem that can be found in Weidmann's paper \cite[Theorem 7]
{Weidmann:Nielsen}, setting $S_{i}=\emptyset $.

\begin{theorem}[Weidmann]
\label{weidman}Let $\Gamma =\left\langle \gamma _{1},\gamma _{2},\ldots
,\gamma _{n}\right\rangle $ be a finitely generated group acting on a tree.
Then either $\Gamma $ is a free group acting freely on $T$ or there is a
Nielsen equivalent set of generators $\Gamma =\left\langle \gamma
_{1}^{\prime },\gamma _{2}^{\prime },\ldots ,\gamma _{n}^{\prime
}\right\rangle $ such that $\gamma _{1}^{\prime }$ is elliptic.
\end{theorem}

\section{The hyperbolic-elliptic case \label{densesec}}

In this section we show that a group generated by an elliptic and a
hyperbolic element is generically dense. In fact, we derive this from the
stronger result, that the group $\Gamma $ generated by a fixed hyperbolic
element and a generic elliptic element is dense. The essence of the proof is
to show that vertex stabilizers of $\Gamma $ are dense in the automorphism
group of the tree rooted at the vertex. This goes by using conjugates of
high powers of the elliptic element by powers of the hyperbolic element.
Although these rooted automorphisms will not be independent random, we shall
be able to pull up sufficient independence to ensure denseness.

\bigskip

Let us start with some notations and basic results on random actions on
rooted trees.

Let $U$ be a spherically homogeneus locally finite infinite rooted tree.
Then $U$ is isomorphic to $T(d_{0},d_{1},\ldots )$ where $d_{i}$ is the
number of children of a vertex at level $i$. For $1\leq l$ let
\begin{equation*}
U_{l}=T(d_{l},d_{l+1},\ldots )
\end{equation*}
denote the subtree of $U$ hanging down from a vertex of level $l$ and let
\begin{equation*}
U^{l}=T(d_{0},d_{1},\ldots ,d_{l-1})
\end{equation*}
be the union of the first $l$ levels of $U$.

Following the elementary decomposition \ref{felbontas}), a Haar uniform
random element $a\in \mathrm{Aut}(U)$ will decompose as
\begin{equation*}
a=\prod_{v\in T}a(v)
\end{equation*}
where $a(v)$ is a uniform random element of the finite group $Sym(Z_{a_{l}})$
with $l=\left| v\right| $. This implies that for all $l\geq 1$ the component
$a_{l}$ is a uniform random element of $\mathrm{Aut}(U^{l})$ and for all
vertices $v\in V$ the component $a_{v}$ is a uniform random element of $%
\mathrm{Aut}(U_{l})$. Also, the set
\begin{equation*}
\left\{ a_{u}\mid u\in U\text{ of level }l\right\}
\end{equation*}
will consist of independent uniform random elements of $\mathrm{Aut}(U_{l})$.

For an element $a\in \mathrm{Aut}(U)$ and a vertex $u\in U$ let $a\circ u$
denote the action of $a^{m}$ on $U_{l}$ where $m$ is the minimal positive
integer such that $a^{m}$ fixes $u$.

Our first lemma already appears in a slightly weaker form in \cite
{AV-dimension_theory}.

\begin{lemma}
\label{unifpro}Let $a$ be a uniform random element of $\mathrm{Aut}(U)$.
Then for each $u\in U$, $a\circ u$ is a uniform random element of $\mathrm{%
Aut}(U_{l})$ where $l=\left| u\right| $.
\end{lemma}

\noindent \textbf{Proof.} Let $u_{0}=u$ and $u_{i}=u_{0}^{a^{i}}$ ($1\leq
i<m $). Let us fix this configuration. Then the $a_{u_{i}}$ ($0\leq i<m$)
are independent random elements of $\mathrm{Aut}(U_{l})$. Now
\begin{eqnarray*}
a^{m} &=&\left( a_{l}\prod_{\left| v\right| =l}a_{v}\right) ^{m}= \\
&=&a_{l}^{m}\left( a_{l}^{-(m-1)}\left( \prod_{\left| v\right|
=l}a_{v}\right) a_{l}^{m-1}\right) \cdots \left( a_{l}^{-1}\left(
\prod_{\left| v\right| =l}a_{v}\right) a_{l}\right) \prod_{\left| v\right|
=l}a_{v}
\end{eqnarray*}
We know that $a_{l}^{m}$fixes $u$ and
\begin{equation*}
a_{l}^{-i}\left( \prod_{\left| v\right| =l}a_{v}\right)
a_{l}^{i}=\prod_{\left| v\right| =l}b_{v}\text{ where }b_{v}=a_{w}\text{
with }w=v^{a_{l}^{i}}
\end{equation*}
This implies that
\begin{equation*}
a\circ u=\left( a^{m}\right) _{u}=\prod_{i=m-1}^{0}a_{u_{i}}
\end{equation*}
the product of independent uniform random elements of $\mathrm{Aut}(U_{l})$,
which is then also uniform random. $\square $

\bigskip

Note that if $u$ and $v$ lie in the same $a$-orbit then $a\circ u$ and $%
a\circ v$ are actually conjugate (being the same product up to a cyclic
permutation), so they are very far from being independent.

We will need to mine out much more independence. To achieve this, we will
take an increasing sequence of vertices $u_{i}$ and look at the action of $%
a\circ u_{i}$ at the first $n$ levels of the tree $U_{l}$.

\begin{lemma}
\label{powerrandom}Let $a$ be a uniform random element of $\mathrm{Aut}(U)$.
Let $u_{1}<u_{2}<\ldots $ be an infinite descending sequence of vertices
such that $\left| u_{i+1}\right| -\left| u_{i}\right| \geq L$ ($i\geq 1$).
Assume that the trees $(U_{u_{i}})^{L}$ of length $L$ are all isomorphic to
\begin{equation*}
V=T(e_{0},e_{1},\ldots ,e_{L-1}).
\end{equation*}
Then the set
\begin{equation*}
\left\{ (a\circ u_{i})_{L}\mid i\geq 1\right\}
\end{equation*}
consists of independent uniform random elements of $\mathrm{Aut}(V)$.
\end{lemma}

\noindent \textbf{Proof.} Let $k_{i}=\left| u_{i}\right| $ ($i\geq 1$). We
have seen in Lemma \ref{unifpro} that $a\circ u_{i}$ is a uniform random
element of $\mathrm{Aut}(U_{u_{i}})$ ($i\geq 1$). This implies that $(a\circ
u_{i})_{L}\in \mathrm{Aut}(V)$ will also be uniform random. Note that $%
(a\circ u_{i})_{L}$ only depends on the values $a(v)$ where $\left| v\right|
<k_{i}+L$, since the elementary automorphisms at higher level vertices fix
the relevant tree $U^{k_{i}+L-1}$.

Let $r\geq 1$. Fix the value of $a_{k_{r}}$. As shown above, this will
determine the value of $(a\circ u_{i})_{L}$ ($i<r$). The proof of Lemma \ref
{unifpro} shows that $(a\circ u_{r})_{L}$ is a function of the set
\begin{equation*}
S=\left\{ a(v)\mid k_{r}\leq \left| v\right| <k_{r}+L\right\} \text{.}
\end{equation*}
Note that the actual function is governed by the orbit structure of $%
a_{k_{r}}$, so it does depend on the values of $a(v)$ ($\left| v\right|
<k_{r}$) but since those are fixed, it is a fixed function of the random
variables in $S$. Now Lemma \ref{unifpro} tells us that $(a\circ
u_{r})_{L}\in \mathrm{Aut}(V)$ is uniform random.

We got that for any fixed value of
\begin{equation*}
(a\circ u_{i})_{L}\text{ (}i<r\text{) }
\end{equation*}
the distribution of $(a\circ u_{r})_{L}$ is uniform random. Using induction
on $r$, this implies that the set
\begin{equation*}
\left\{ (a\circ u_{i})_{L}\mid 1\leq i\leq r\right\}
\end{equation*}
consists of independent uniform random elements of $\mathrm{Aut}(V)$. $%
\square $

\bigskip

We need the following technical lemma on rooted trees.

\begin{lemma}
\label{techno}Let $T=T(e_{0},e_{1},\ldots ,e_{L-1})$ be a rooted tree with $%
e_{0}\geq 3$. Let $x$ and $y$ be distinct vertices of $T$ on the first
level. Let
\begin{equation*}
X=\left\{ g\in \mathrm{Aut}(V)\mid zg=z\text{ }\forall z\geq x\right\}
\end{equation*}
and
\begin{equation*}
Y=\left\{ g\in \mathrm{Aut}(V)\mid zg=z\text{ }\forall z\geq y\right\}
\end{equation*}
Then $X$ and $Y$ generate $\mathrm{Aut}(V)$.
\end{lemma}

\noindent \textbf{Proof.} Let $v\in T$ be a vertex distinct from the root.
Then $v\ngeq x$ or $v\ngeq y$ so all elementary automorphisms at $v$ are
contained in $X\cup Y$. Since $e_{0}\geq 3$, the symmetric group $\mathrm{Sym%
}(e_{0})$ is generated by any two distinct point stabilizers. This implies
that any elementary automorphism at the root is generated by elements of $X$
and $Y$. Thus all elementary automorphisms are generated by $X\cup Y$ and
the Lemma holds. $\square $

\bigskip

Now we start to discuss random generation for unrooted trees. For a vertex $%
t\in T$ let $\mathrm{Aut}(T)_{t}$ denote the stabilizer of $t$ in $\mathrm{%
Aut}(T)$.

\begin{figure}[tbp]
\includegraphics[width=12cm]{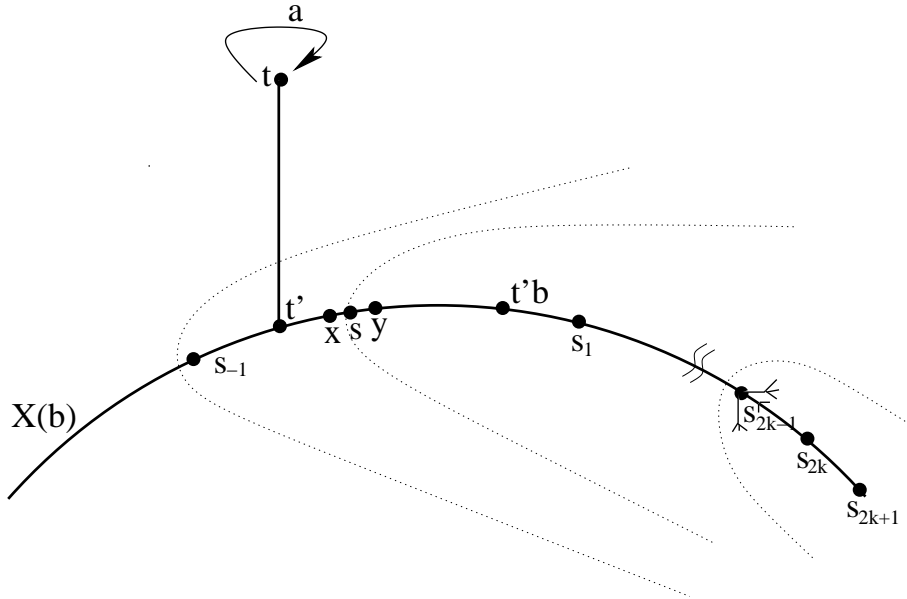}\newline
\caption{Mining out randomness}
\label{fig:mining}
\end{figure}

\begin{lemma}
\label{densepoint}Let $b\in \mathrm{Aut}(T)$ be hyperbolic and let $t\in T$
be a vertex. Let $a\in \mathrm{Aut}(T)_{t}$ be a random element and let
\begin{equation*}
\Gamma =\left\langle a,b\right\rangle \leq \mathrm{Aut}(T)\text{.}
\end{equation*}
Then generically, there exists $s\in X(b)$ of degree at least $3$ such that
the stabilizer $\Gamma _{s}$ is dense in $\mathrm{Aut}(T)_{s}$.
\end{lemma}

\noindent \textbf{Proof.} Let $l=l(b)$ be the translation distance of $b$.
Let $t^{\prime }$ be the projection of $t$ on $X(b)\ $and let $s$ be an
element of the geodesic $(t^{\prime },t^{\prime }b]$ with degree at least $3$%
. Let $x\in \lbrack t^{\prime },s]$ with $d(x,s)=1$ and let $y\in \lbrack
s,sb]$ with $d(s,y)=1$. Let $T_{s}$ be the tree $T$ rooted at $s$.

Let us fix a positive integer $K$. We claim that generically, the action of $%
\Gamma _{s}$ on $T_{s}^{K}$ equals the full automorphism group $\mathrm{Aut}%
(T_{s}^{K})$. For an integer $n$ let
\begin{equation*}
s_{n}=sb^{Kn},\text{ let }x_{n}=xb^{Kn}\text{ and let }y_{n}=yb^{Kn}
\end{equation*}
and let
\begin{equation*}
V_{i}=\mathrm{Shadow}_{x_{2i-1}\longrightarrow s_{2i-1}}\text{ and }W_{i}=%
\mathrm{Shadow}_{y_{2i-1}\longrightarrow s_{2i-1}}
\end{equation*}
as trees rooted at $s_{2i-1}$.

Let us apply Lemma \ref{powerrandom} for the tree $U=T$ (rooted at $t$),
setting $L=2K$ and $u_{i}=s_{2i-1}$ ($i\geq 1$). We get that the set
\begin{equation*}
\left\{ (a\circ u_{i})_{L}\mid i\geq 1\right\}
\end{equation*}
consists of independent uniform random elements of $\mathrm{Aut}(V_{i}^{L})$.

Let $g$ be an arbitrary element of $\mathrm{Aut}(V_{-1}^{L})$ that fixes
every vertex of $V_{-1}^{L}$ that is not a descendant of $s$. Then because
of the above independence, generically there exists $i_{g}\geq 1$ such that
\begin{equation*}
a_{g}=b^{2Ki_{g}}(a\circ u_{i_{g}})b^{-2Ki_{g}}\in \Gamma _{s}
\end{equation*}
acts on $V_{-1}^{L}$ as $g$ does. Let $z\in T_{s}^{K}$ with $z\geq x$. Then $%
d(z,s)\leq K$ and $d(s_{-1},s)=K$ so $d(s_{-1},z)\leq 2K$, which implies $%
za_{g}=z$.

The same way, applying Lemma \ref{powerrandom} for $u_{i}=s_{-2i+1}$ ($i\geq
1$) we get that for an arbitrary element $h\in \mathrm{Aut}(W_{1}^{L})$ that
fixes every vertex of $W_{1}^{L}$ that is not a descendant of $s$,
generically there exists $i_{h}\geq 1$ such that
\begin{equation*}
a_{h}=b^{-2Ki_{h}}(a\circ u_{i_{h}})b^{2Ki_{h}}\in \Gamma _{s}
\end{equation*}
acts on $W_{1}^{L}$ as $h$ does. Also, for all $z\in T_{s}^{K}$ with $z\geq
y $ we have $za_{h}=z$.

Applying Lemma \ref{techno} on the tree $T_{s}^{K}$ we get that the actions
of all the $a_{g}$ and $a_{h}$ on $T_{s}^{K}$ generate $\mathrm{Aut}%
(T_{s}^{K})$ and so our claim holds.

Since for all $K$ the claim generically holds, it generatically holds for
all $K$ at the same time. In particular, $\Gamma _{s}$ is generically dense
in $\mathrm{Aut}(T)_{s}$. $\square $

\bigskip

The following lemma will be used to establish densenesss in $\mathrm{Aut}%
^{0}(T)$. To make the notation simpler, we will use regular trees here (that
is, we forget the baricentric points).

\begin{figure}[tbp]
\includegraphics[width=10cm]{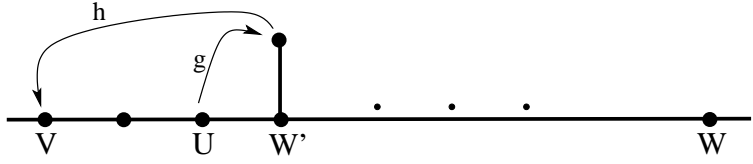}\newline
\caption{Transitive action on the even vertices}
\label{fig:tr_action}
\end{figure}

\begin{lemma}
\label{kutykuty}Let $\Gamma $ be a subgroup of $\mathrm{Aut}(T)$. Assume
that $\Gamma $ contains a hyperbolic element and that there exists $s\in T$
such that $\Gamma _{s}$ is dense in $\mathrm{Aut}(T)_{s}$. Then the closure
of $\Gamma $ contains $\mathrm{Aut}^{0}(T)$.
\end{lemma}

\noindent \textbf{Proof.} Let $b$ be a hyperbolic element of $\Gamma $. Let $%
X$ be the set of vertices that have even distance from $s$. Let $Y$ be the $%
\Gamma $-orbit of $s$. Using conjugation, it is easy to see that for all $%
y\in Y$, the stabilizer $\Gamma _{y}$ is dense in $\mathrm{Aut}(T)_{y}$.

We claim that $X\subseteq Y$. To show this, let $u\in X\cap Y$ and let $v\in
T$ such that $d(u,v)=2$. Let $w=s$ or $w=sb$ such that $w\neq u$. Let $%
w^{\prime }\in \lbrack u,w]$ adjacent to $u$ and let $v^{\prime }\notin
\lbrack u,w]$ be a vertex that is adjacent to $w^{\prime }$. Such $v^{\prime
}$ exists since the degree of $w^{\prime }$ is at least $3$. Now using
density at $w$ and $d(w,u)=d(w,v^{\prime })$, there exists an element $g\in
\Gamma _{w}$ with $ug\in v^{\prime }$. This implies $v^{\prime }\in X\cap Y$%
. Now $d(u,v)=d(u,v^{\prime })=2$ so using density again at $u$, there
exists $h\in \Gamma _{u}$ with $v^{\prime }h\in v$, implying $v\in X\cap Y$.
Thus every element of distance $2$ from an element in $X\cap Y$ lies in $Y$
and the claim holds.

Let $G$ be the closure of $\Gamma $ in $\mathrm{Aut}(T)$. The group $\mathrm{%
Aut}^{0}(T)$ acts transitively on $X$ and for all $x\in X$ the stabilizers
\begin{equation*}
G_{x}=\mathrm{Aut}^{0}(T)_{x}=\mathrm{Aut}(T)_{x}
\end{equation*}
So if $X=Y$ then $G=\mathrm{Aut}^{0}(T)$, otherwise $G=\mathrm{Aut}(T)$. $%
\square $

\bigskip

We are ready to establish the main result of this section. Note that we do
use the baricentric subdivision here, that is, we allow vertices of degree $%
2 $.

\begin{theorem}
\label{ellhyp}Let $a$ be a generic elliptic and $b\mathrm{\ }$be a generic
hyperbolic element of $\mathrm{Aut}(T)$. Then
\begin{equation*}
\mathrm{Aut}^{0}(T)\subseteq \overline{\langle a,b\rangle }.
\end{equation*}
\end{theorem}

\noindent \textbf{Proof.} Let $b$ be a fixed hyperbolic element of $\mathrm{%
Aut}(T)$, let $t\in T$ and let $a$ be a generic element of $\mathrm{Aut}%
(T)_{t}$. Let $\Gamma =\langle a,b\rangle $ and let $G$ be the closure of $%
\Gamma $ in $\mathrm{Aut}(T)$. Using Lemma \ref{densepoint} there exists $%
s\in X(b)$ of degree at least $3$ such that the stabilizer $\Gamma _{s}$ is
dense in $\mathrm{Aut}(T)_{s}$. Since $s$ is not a baricenter, we can use
Lemma \ref{kutykuty} which implies $\mathrm{Aut}^{0}(T)\subseteq G$.

The measure theoretic part of the theorem now follows by summing over $t$
and integrating over the hyperbolic coordinate. Let $\mathrm{Ell}$ denote
the set of elliptic and $\mathrm{Hyp}$ the set of hyperbolic elements in $%
\mathrm{Aut}(T)$. By Baire's theorem, to establish the topological
statement, it is enough to show that for any $\phi \in \mathrm{Aut}^{0}(T)$
and any $n\in N$ the set
\begin{equation*}
W(n,\phi )=\left\{ (a,b)\in \mathrm{Hyp}\times \mathrm{Ell}\left| \exists
\psi \in {\langle a,b\rangle }\text{ such that }\phi |B(t_{0},n)=\psi
|B(t_{0},n)\right. \right\}
\end{equation*}
is a dense open subset of the set $\mathrm{Hyp}\times \mathrm{Ell}$. Recall
that $B(t,n)$ denotes the ball of radius $n$ around $t$. The set $W(n,\phi )$
is open, because if the word $\psi =w(a,b)$ realizes the condition $\phi
|B(t_{0},n)=\psi |B(t_{0},n)$ then so will $\psi ^{\prime }=w(a^{\prime
},b^{\prime })$ if $a^{\prime }$ and $b^{\prime }$ are sufficiently close to
$a$ and $b$, respectively. The density of $W(n,\phi )$ follows from the
measure theoretic part.

The theorem holds. $\square $

\section{Trichotomy \label{trichotomysec}}

In this section we establish Theorem \ref{thm:main}. First we need some
notations.

\bigskip

Let
\begin{eqnarray*}
\mathrm{Ell}_{n} &=&\left\{ \mathbf{a}\in A^{n}\mid \mathbf{a}_{i}\text{ is
elliptic }(1\leq i\leq n)\right\} \\
\mathrm{Hyp}_{n} &=&\left\{ \mathbf{a}\in A^{n}\mid \mathbf{a}_{i}\text{ is
hyperbolic }(1\leq i\leq n)\right\} \\
\mathrm{Mix}_{n} &=&A^{n}\backslash (\mathrm{Ell}_{n}\cup \mathrm{Hyp}_{n})
\end{eqnarray*}

Also let
\begin{equation*}
C=\left\{ \mathbf{a}\in \mathrm{Ell}_{n}\mid \left\langle \mathbf{a}%
\right\rangle \text{ is precompact}\right\}
\end{equation*}
\begin{equation*}
S=\left\{ \mathbf{a}\in \mathrm{Hyp}_{n}\mid \mathbf{a}\text{ satisfies the
Schottky-condition}\right\}
\end{equation*}
\begin{equation*}
D=\left\{ \mathbf{a}\in \mathrm{Mix}_{n}\mid \left\langle \mathbf{a}%
\right\rangle \text{ is dense in }A\text{ or }A^{0}\right\}
\end{equation*}

It is easy to see that $C,S$ and $D$ are pairwise disjoint, $C$ is open and
closed and $S$ is open.

\begin{lemma}
\label{kicsi}$D$ is generic and co-meager in $\mathrm{Mix}_{n}$\textrm{,}
that is, $\mathrm{Mix}_{n}\backslash D$ is meager of measure $0$.
\end{lemma}

\noindent \textbf{Proof.} For $n=2$ this is Theorem \ref{ellhyp}. For $n>2$
and $1\leq i,j\leq n$ let
\begin{equation*}
E_{ij}=\left\{ \mathbf{a}\in \mathrm{Mix}_{n}\backslash D\mid \mathbf{a}_{i}%
\text{ is hyperbolic and }\mathbf{a}_{j}\text{ is elliptic}\right\}
\end{equation*}
Using Theorem \ref{ellhyp} the projection of $E_{ij}$ to the coordinates $%
i,j $ is meager of measure zero in $A^{2}$. So $E_{ij}$ is also meager of
measure zero. The lemma now follows from
\begin{equation*}
\mathrm{Mix}_{n}\backslash D=\bigcup_{1\leq i,j\leq n}E_{ij}\text{.}
\end{equation*}
$\square $

\bigskip

Let $\Delta $ denote the Nielsen action of $\mathrm{Aut}(F_{n})$ on $A^{n}$.
Since the action does not change the subgroup generated by the tuple, $C$
and $D$ are invariant under $\Delta $.

\begin{lemma}
We have $S^{\Delta }\subseteq \mathrm{Hyp}_{n}$.
\end{lemma}

\noindent \textbf{Proof.} If $\mathbf{a}$ satisfies the Schottky-condition
then $\left\langle \mathbf{a}\right\rangle $ is discrete and free. That is,
vertex stabilizers are finite, which, using that the free group is
torsion-free, implies that they are trivial. So $\left\langle \mathbf{a}%
\right\rangle $ contains no elliptic elements. The same holds for $%
\left\langle \mathbf{a}^{\varphi }\right\rangle =\left\langle \mathbf{a}%
\right\rangle $, in particular, all the entries of $\mathbf{a}^{\varphi }$
are hyperbolic. $\square $

\bigskip

\begin{lemma}
\label{ize}We have
\begin{equation*}
S^{\Delta }=\left\{ \mathbf{a}\in A^{n}\mid \left\langle \mathbf{a}%
\right\rangle \text{ is discrete and free}\right\} .
\end{equation*}
\end{lemma}

\noindent \textbf{Proof.} This follows from Lubotzky's theorem on the
Schottky condition (see Section 2) and the fact that any two minimal
generating sets of a free group are Nielsen-equivalent. $\square $

\bigskip

\begin{lemma}
\label{egylepes}For $\mathbf{a}\in \mathrm{Ell}_{n}\backslash C$ there
exists a Nielsen transformation $\varphi \in \Delta $ such that $\mathbf{a}%
^{\varphi }\in \mathrm{Mix}_{n}$.
\end{lemma}

\noindent \textbf{Proof.} Since $\mathbf{a}\notin C$, there is no common
fixed points for the $\mathbf{a}_{i}$. Since the fixed-point sets of the $%
\mathbf{a}_{i}$ are convex, it follows from the Caratheodory theorem on
trees that there is $i,j$ such that $\mathbf{a}_{i}$ and $\mathbf{a}_{j}$
have no common fixed points. But this yields that $\mathbf{a}_{i}\mathbf{a}%
_{j}$ is hyperbolic, implying $\mathbf{a}^{R_{i,j}^{+}}\in \mathrm{Mix}_{n}$%
. $\square $

\bigskip

\begin{lemma}
\label{soklepes}For $\mathbf{a}\in \mathrm{Hyp}_{n}\backslash S^{\Delta }$
there exists $\varphi \in \Delta $ such that $\mathbf{a}^{\varphi }\in
\mathrm{Mix}_{n}$.
\end{lemma}

\noindent \textbf{Proof.} Lemma \ref{ize} implies that $\left\langle \mathbf{%
a}\right\rangle $ is not discrete and free. Then the theorem of Weidmann
(Theorem \ref{weidman}) tells us that there is $\xi \in \Delta $ such that
the tuple $\mathbf{a}^{\xi }$ contains an elliptic element. If $\mathbf{a}%
^{\xi }\in \mathrm{Mix}_{n}$, we proved our lemma. If $\mathbf{a}^{\xi }\in
\mathrm{Ell}_{n}$ then $\mathbf{a}^{\xi }\notin C$, otherwise $\left\langle
\mathbf{a}^{\xi }\right\rangle $ would consist of elliptic elements,
contradicting $\mathbf{a}\in \mathrm{Hyp}_{n}$. So $\mathbf{a}^{\xi }\in
\mathrm{Ell}_{n}\backslash C$ and using Lemma \ref{egylepes} we see that
there exists $\delta \in \Delta $ such that $\mathbf{a}^{\xi \delta }\in
\mathrm{Mix}_{n}$. $\square $

\bigskip

We are ready to prove Theorem 1.

\bigskip

\noindent \textbf{Proof of Theorem 1.} Let $\mathbf{a}\in A^{n}$ be a
generic $n$-tuple. First we show that $\mathbf{a}$ generates a free group of
rank $n$. The measure generic part directly follows from a result of the
first author \cite[Corollary 1.6]{abert}. For the topological part, let $w$
be a nontrivial word in $n$ letters. Then the support
\begin{equation*}
\mathrm{Supp}(w)=\left\{ \mathbf{a}\in A^{n}\mid w(a)=1\right\}
\end{equation*}
is closed in $A^{n}$ and by the measure theoretic part, it has measure zero.
This implies that it is nowhere dense. Hence the set of points satisfying
any nontrivial words is meager and so the topological part follows.

Now we establish the trichotomy. Let
\begin{equation*}
L=(\mathrm{Mix}_{n}\backslash D)^{\Delta }.
\end{equation*}
Then using Lemma \ref{kicsi}, $L$ is a union of countably many $\mathrm{Mix}%
_{n}\backslash D$-translates, so $\mu (L)=0$ and $L$ is meager.

Let $\mathbf{a\notin }L$. We claim that there exists $\varphi \in \Delta $
such that $\mathbf{a}^{\varphi }\in C\cup D\cup S$.

If $\mathbf{a}\in \mathrm{Ell}_{n}$ and $\mathbf{a}\in C$ we are done. If $%
\mathbf{a}\in \mathrm{Ell}_{n}\backslash C$ then by Lemma \ref{egylepes}
there exists a Nielsen transformation $\varphi \in \Delta $ such that $%
\mathbf{a}^{\varphi }\in \mathrm{Mix}_{n}$. However, $\mathbf{a}^{\varphi
}\notin L\supseteq \mathrm{Mix}_{n}\backslash D$ implying $\mathbf{a}%
^{\varphi }\in D$.

If $\mathbf{a}\in \mathrm{Hyp}_{n}$ and $\mathbf{a}\in S^{\Delta }$ then
there exists a $\varphi \in \Delta $ such that $\mathbf{a}^{\varphi }\in S$.
Otherwise $\mathbf{a}\in \mathrm{Hyp}_{n}\backslash S^{\Delta }$ and by
Lemma \ref{soklepes} there exists $\varphi \in \Delta $ such that $\mathbf{a}%
^{\varphi }\in \mathrm{Mix}_{n}$. Again, $\mathbf{a}^{\varphi }\notin
L\supseteq \mathrm{Mix}_{n}\backslash D$ so $\mathbf{a}^{\varphi }\in D$.

Finally, if $\mathbf{a}\in \mathrm{Mix}_{n}$ then $\mathbf{a}\notin
L\supseteq \mathrm{Mix}_{n}\backslash D$ so $\mathbf{a}\in D$. We have
proved our claim.

Since the group generated by $\mathbf{a}^{\varphi }$ equals the group
generated by $\mathbf{a}$, the trichotomy holds. $\square $

\bigskip

\noindent \textbf{Remark. }The freeness of a generic subgroup implies that
the action of $\mathrm{Aut}(F_{n})$ on $\mathrm{Aut}(T)^{n}$ is essentially
free. Another consequence is that apart from a nullset in $\mathrm{Aut}%
(T)^{n}$, the $\mathrm{Aut}(F_{n})$-orbits can be identified with the
subgroup that an element of the orbit generates. This allows us to talk
about \emph{generic subgroups} rather than the group generated by a generic $%
n$-tuple.

\section{Applications \label{primitivsec}}

In this section we use the existence of dense free subgroups in $\mathrm{Aut}%
^{0}(T)$ to find some interesting permutation actions on countable sets.
Then we discuss why the trichotomy theorem does not hold in general in the
realm of locally finite groups.

Let us recall some notions of group actions. Let the group $\Gamma $ act on
the set $\Omega $. The action is $k$\emph{-transitive}, if the action of $%
\Gamma $ on the ordered $k$-tuples of distinct elements of $\Omega $ is
transitive. The action is \emph{primitive} if there is no non-trivial $%
\Gamma $-invariant equivalence relation on $\Omega $. This is equivalent to
say that the action is transitive and a point stabilizer $\Gamma _{\omega }$
is a maximal subgroup in $\Gamma $. The action of $\Gamma $ is \emph{%
quasi-primitive} if every normal subgroup of $\Gamma $ is acts transitively
or trivially. It is easy to see that every $2$-transitive action is
primitive and every primitive action is quasi-primitive. Primitive actions
of finite and infinite groups have a well-established theory.
Quasi-primitive actions of finite groups have also been extensively studied
(see \cite{preeger} and references therein).

Let $G$ be a group. A subgroup $H\leq G$ is \emph{subnormal}, if there is a
chain
\begin{equation*}
H=G_{n}\lhd G_{n-1}\lhd \ldots \lhd G_{0}=G.
\end{equation*}
We say that the action of $\Gamma $ is \emph{subnormal transitive} if every
nontrivial subnormal subgroup of $\Gamma $ acts transitively or trivially.
It is easy to see that primitive actions are not necessarily subnormally
transitive, but if we assume that the group is $k$-transitive for every $k$
then it is also subnormal transitive.

The following general lemma will establish subnormal transitivity for a
general class of actions.

\begin{lemma}
\label{subnorm}Let $G$ be a disconnected topologically simple topological
group and let $O<G$ be an open subgroup. Let $\Gamma <G$ be a dense subgroup
and let $\Lambda =\Gamma \cap O$. Then the right coset action of $\Gamma $
on $\Gamma /\Lambda $ is faithful and subnormal transitive.
\end{lemma}

\noindent \textbf{Proof.} Let $N\neq 1$ be a subnormal subgroup of $\Gamma $%
. Then the closure of $N$ in $G$ is subnormal in the closure of $\Gamma $ in
$G$, which is equal to $G$, since $\Gamma $ is dense in $G$. But $G$ is
topologically simple, implying that $N$ is dense in $G$. Now using that $O$
is open, we get $NO=G$. This implies
\begin{equation*}
N\Lambda =N(\Gamma \cap O)=\Gamma
\end{equation*}
which is equivalent to saying that $N$ acts transitively on $\Gamma /\Lambda
$. $\square $

\bigskip

\noindent \textbf{Proof of Corollary \ref{szep}.} Let $T$ be the $3$-regular
infinite tree. Let $\Gamma $ be a generic $n$-generated subgroup of $\mathrm{%
Aut}(T)$. Then using Theorem \ref{thm:main}, $\Gamma $ is free on $n$
generators and on a set of infinite measure, its closure in $\mathrm{Aut}(T)$
equals $\mathrm{Aut}^{0}(T)$. In particular, there exists a subgroup $\Gamma
\subseteq \mathrm{Aut}^{0}(T)$ which is isomorphic to $F_{n}$ and dense in $%
\mathrm{Aut}^{0}(T)$.

Let $t\in T$ be a vertex, let $G=\mathrm{Aut}^{0}(T)$ and let $O=\mathrm{Aut}%
^{0}(T)_{t}$. Applying Lemma \ref{subnorm} we get that the right coset
action of $\Gamma $ on $\Gamma /\Lambda $ is faithful and subnormal
transitive, where $\Lambda =\Gamma _{t}$ is the vertex stabilizer of $\Gamma
$. This right coset action is permutation isomorphic to the action of $%
\Gamma $ on $T$ and $\Gamma $ is dense in $\mathrm{Aut}^{0}(T)$, so the
action of $\Gamma $ is primitive but not $2$-transitive. For every finite
subset $X\subseteq T$, the pointwise stabilizer of $X$ in $\mathrm{Aut}%
^{0}(T)$ is nontrivial, hence using density again, the same holds for $%
\Gamma $. The corollary holds. $\square $

\bigskip

Now we show that the trichotomy theorem fails for the product of two trees.

\begin{proposition}
\label{bzzz}Let $T$ and $U$ be regular trees, both with degree at least $3$.
Let $A=\mathrm{Aut}(T)\times \mathrm{Aut}(U)$. Then for every $n\geq 2$
there is a subset $X\subseteq A^{n}$ of infinite measure such that for all $%
\mathbf{a}\in X$, the closure of the subgroup generated by the elements of $%
\mathbf{a}$ is not discrete, open or compact.
\end{proposition}

\noindent \textbf{Proof.} Let $\pi _{T}:A\rightarrow \mathrm{Aut}(T)$ and $%
\pi _{U}:A\rightarrow \mathrm{Aut}(U)$ be the projections to the two
coordinates. Let
\begin{equation*}
X=\left\{ \mathbf{a}\in A^{n}\mid \pi _{T}(\left\langle \mathbf{a}%
\right\rangle )\text{ is dense in }\mathrm{Aut}(T)\text{ and }\overline{\pi
_{U}(\left\langle \mathbf{a}\right\rangle )}\text{ is compact}\right\}
\end{equation*}
The normalized Haar measure on $A$ equals the product of the normalized Haar
measures on $\mathrm{Aut}(T)\times \mathrm{Aut}(U)$, so Theorem \ref
{thm:main} implies that $X$ has infinite measure.

Let $\mathbf{a}\in X$, let $\Gamma =\left\langle \mathbf{a}\right\rangle $
and let $G$ be the closure of $\Gamma $ in $A$. We claim that $G$ does not
satisfy the trichotomy. Indeed, $G$ can not be compact, since the continuous
image $\pi _{T}(G)\supseteq \pi _{T}(\Gamma )$ is dense in $\mathrm{Aut}(T)$
and hence is not compact. Using density again, there exists a sequence $%
(\gamma _{i})$ of distinct elements of $\Gamma $ such that $(\pi _{T}(\gamma
_{i}))$ is convergent. Since $\overline{\pi _{U}(\Gamma )}$ is compact,
there is a subsequence $(\delta _{i})$ of $(\gamma _{i})$ such that $(\pi
_{U}(\delta _{i}))$ is convergent. But then $(\delta _{i})$ is convergent,
so $\Gamma $ (and hence $G$) is not discrete.

Finally, assume by contradiction that $G$ is open. Then $O=\pi _{U}(G)$ is
an open subgroup of $\mathrm{Aut}(U)$ and since $\Gamma $ is dense in $G$, $%
\pi _{U}(\Gamma )$ is dense in $O$. This implies that $O$ is topologically
finitely generated. The group $O$ is also compact, being an open subgroup of
$\overline{\pi _{U}(\Gamma )}$. But then $O$ must stabilize a vertex or a
geometric edge. This means that there is an open subgroup of $O$ of index at
most $2$ that fixes a vertex $t\in T$, which implies that the vertex
stabilizer $\mathrm{Aut}(U)_{t}$ has a topologically finitely generated open
subgroup. Since $\mathrm{Aut}(U)_{t}$ is compact, we get that $\mathrm{Aut}%
(U)_{t}$ is itself topologically finitely generated, a contradiction. Hence $%
G$ is not open. $\square $

\end{document}